\documentstyle[12pt]{article}
\newtheorem{definition}{Definition}

\newtheorem{theorem}{Theorem}
\newtheorem{proposition}{Proposition}
\newtheorem{corollary}{Corollary}
\newtheorem{property}{Property}

\begin{document}

\author{Marc Mars$^\natural$\thanks{E-mail address: M.Mars@qmw.ac.uk} 
\ and Roustam M. Zalaletdinov$^{\flat \sharp}$\thanks{E-mail address: 
zala@camk.edu.pl} \thanks{\'Sniadeccy Fellow} 
\thanks{Address after 10 May 1997:
School of Mathematical Sciences, Queen Mary \& Westfield College, 
University of London, Mile End Road, London E1 4NS, England, U.K. 
E-mail address: rmz@maths.qmw.ac.uk} \\[5mm]
{\em $^\natural$School of Mathematical Sciences, Queen Mary \& Westfield College}\\
{\em University of London, Mile End Road, London E1 4NS, England, U.K.}\\[2mm]
{\em $^\flat$N. Copernicus Astronomical Center, Polish Academy of Sciences}\\
{\em  ul. Bartycka 18, 00-716 Warsaw, Poland}\\[2mm]
{\em $^\sharp$Department of Theoretical Physics, Institute of Nuclear Physics}\\
{\em Uzbek Academy of Sciences, Tashkent 702132, Uzbekistan, C.I.S.}}
\title{\LARGE \bf Space-time Averages in Macroscopic Gravity \\
and Volume-preserving Coordinates}
\date{31 January 1997}

\maketitle

\begin{abstract}
The definition of the covariant space-time averaging scheme for the objects
(tensors, geometric objects, etc.) on differentiable metric manifolds with a 
volume $n$-form, which has been 
proposed for the formulation of macroscopic gravity, is analyzed. 
An overview of the space-time averaging procedure in Minkowski spacetime
is given and comparison between this averaging scheme and that 
adopted in macroscopic gravity is carried out throughout the paper.
Some
new results concerning the algebraic structure of the   
averaging operator are precisely formulated and proved, the main one being
that the averaging bilocal operator is idempotent iff it is factorized
into a bilocal product of a matrix-valued function on the manifold, taken 
at a
point, by its inverse at another point. The previously proved existence
theorems for the averaging and coordination bilocal operators are revisited 
with more detailed proofs of related results. A number of new results
concerning the structure of the volume-preserving averaging operators and 
the class of proper coordinate systems are
given. It is shown, in particular, that such operators are defined on an
arbitrary $n$-dimensional differentiable metric manifold 
with a 
volume $n$-form up to the freedom of $
(n-1) $ arbitrary functions of $n$ arguments and 1 arbitrary function of $
(n-1)$ arguments. All the results given in this paper are also valid
whenever appropriate for affine connection manifolds including
(pseudo)-Riemannian manifolds.
\end{abstract}
\thispagestyle{empty}

\newpage

\section{Introduction}

Space-time averaging procedures play an important role in modern 
physics because of, first of all, 
their relevance to deriving classical macroscopic
theories. A well-known example of such a procedure is the space-time scheme
developed for averaging out the microscopic Lorentz electrodynamics to
derive the macroscopic Maxwell electrodynamics (see, for example, \cite
{Nova:1955}-\cite{Inga-Jami:1985})\footnote{
Alternative approaches for deriving macroscopic electrodynamics
which apply other averaging procedures (for example, space averaging and
statistical ensemble averaging) are not considered here. For such approaches
and a discussion about their interrelations, physical significance, etc. see 
\cite{deGr:1969}-\cite{Zala:1996}.}. Another important physical argument for
considering space-time averaging procedures is that 
they are relevant to modelling the process of physical measurement. It
is the space-time averages of physical fields that are known \cite
{Bohr-Rose:1933},\cite{DeWi:1962} to have direct observational status and
physical meaning.

The space-time averaging procedure of classical electrodynamics utilizes
explicitly the flat character of the Minkowski space-time manifold and its
formulation is essentially based on the existence of Cartesian coordinates.
In this connection, the following questions are of interest and importance:
whether it is possible ({\em a}) to encode the properties of the averaging
procedure in a covariant manner suitable for differentiable manifolds not
equipped
with a metric and a connection in general, and ({\em b}) to formulate the
corresponding covariant averaging procedure for the objects which can live
on such manifolds and to clear up its geometric meaning. Having replied
({\em a}) and ({\em b}), the next question of primary importance is ({\em 
c}) how the averaging procedure may be made compatible with a metric and a
connection when the differentiable manifold possesses those structures.

It should be pointed out here that the above problem of generalizing 
the flat space-time procedure for curved manifolds goes far beyond being
simply an
academic problem. One of its most important areas of applicability
is the general theory of relativity where
the space-time is a 4-dimensional pseudo-Riemannian manifold. As
is well-known,
there is not yet a satisfactory derivation of the Maxwell equations in
general relativity, apart from the covariantization procedure where the
partial derivatives of the special theory of relativity are replaced by
covariant ones
(see, for example, \cite
{MTW:1973}). A physically motivated and mathematically correct derivation of
the macroscopic Maxwell equations by averaging out the
general relativistic microscopic Maxwell-Lorentz equations is still
lacking and, furthermore,
the very foundations of microscopic electrodynamics in general
relativity are not well-established yet. Another problem in general relativity
where the availability of a space-time averaging procedure is of primary
importance is the so-called averaging problem (see \cite{Zala:1996},\cite
{Elli:1984}-\cite{Kras:1996} for a review and discussion). Its main
motivation comes from cosmology where Einstein's equations
are usually utilized with a hydrodynamic stress-energy tensor without any
satisfactory proof of why the left-hand side of the equations (the field
operator) keeps the same structure while the right-hand side has been
changed, or averaged, from a real discrete matter distribution
(stars, galaxies, etc.) to a continuous one. The task here is to carry out a
space-time averaging of Einstein's equations in order to
understand the structure of the averaged (macroscopic) field equations and 
apply them to deal with the overwhelming majority of cosmological
problems. A solution for this problem is also desirable 
in order to provide a rigorous basis for
constructing continuous matter models from discrete ones
(this construction is based in modern cosmology mainly on phenomenological
grounds).

The goal of this paper is twofold: (1) to give precise formulations and 
detailed proofs
of some new
results concerning the algebraic structure of the averaging and
coordination 
operators and the
properties of the space-time averages within the so-called
macroscopic gravity \cite{Zala:1996},\cite{Zala:1992},\cite{Zala:1993}-\cite
{Zala:press} (the basics of this averaging procedure has been developed in \cite
{AZS:1990a},\cite{AZS:1990b}); (2) to revisit the previously
shown existence of the bilocal operators of macroscopic gravity \cite
{Zala:1992},\cite{Zala:1993} with more detailed proofs of related results.
The paper is organized as follows. Section 2 gives an overview of the
space-time averaging scheme used in electrodynamics. The space-time averages
of macroscopic gravity are defined in Section 3. The next two sections are
devoted to the algebra of the averaging operator
with a discussion of its algebraic properties and
a formulation of an important theorem stating
that the averaging bilocal operator is
idempotent iff it is factorized into a bilocal product of a matrix-valued
function on the manifold, taken at a point, by its inverse at another
point. The differential properties of the averages and the coordination 
bivector are summarized in Section
6. In Section 7 the existence theorems for the averaging and 
coordination operators are formulated and proved. The last Section 8 
is devoted to the definition and properties of the
proper coordinate systems, the coordinates in which the
averaging operators take their simplest form. It is shown, in particular,
that such operators are defined on an arbitrary $n$-dimensional
differentiable metric manifold with a volume $n$-form up to the 
freedom of $(n-1)$ arbitrary functions of $
n$ arguments and one arbitrary function of $(n-1)$ arguments. 
A summary on the meaning and content of the space-time averaging procedure of 
macroscopic gravity is given in the Conclusion. 

Though given for $n$-dimensional differentiable metric manifolds with a 
volume $n$-form, 
they are valid whenever
appropriate for affine connection manifolds, including 
(pseudo)-Riemannian manifolds.
All the results in this paper also hold generally 
for
$n$-dimensional differentiable manifolds with a 
volume $n$-form possessing neither a metric nor a connection. 

\section{Space-time averages in Minkowski manifold}

The space-time averaging procedure applied in electrodynamics is explicitly
based on the flat character of the Minkowski space-time ${\cal E}$ and the
existence of the Cartesian coordinates $(t,x^a),a=1,2,3$. The
space-time averages are defined as follows \cite{Nova:1955}-\cite
{Inga-Jami:1985}.

\begin{definition}
\label{average:ED} The average value of a tensor field $p_\beta ^\alpha
(x),\,x\in {\cal E}$, over a compact space region $S$ and a finite time
interval $\Delta t$ at a supporting point $(t,x^a)\in \Delta t\times S$ is 
\begin{equation}
\label{defaver:ED}\langle p_\beta ^\alpha (t,x^a)\rangle _{{\cal E}}=\frac
1{\Delta tV_S}\int_{\Delta t}\int_Sp_\beta ^\alpha (t+t^{\prime
},x^a+x^{a\prime })dt^{\prime }d^3x^{\prime }\ .
\end{equation}
Here $V_S$ is the 3-volume of the region $S$, which is usually
taken as a
3-sphere of radius $R$ around the point $x^a$ at the instant of time $t$, 
$$
V_S=\int_Sd^3x^{\prime }\ . 
$$
\end{definition}

Formula (\ref{defaver:ED}) defines the average value $\langle p_\beta
^\alpha (t,x^a)\rangle _{{\cal E}}$ at a point $(t,x^a)$ on the manifold $
{\cal E}$. In order to obtain averaged tensor fields and define its 
derivatives, it is necessary to make additional assumptions
concerning the averaging regions $S$ and the intervals $\Delta t$. This
assumptions are usually made only tacitly (see, however, a discussion in \cite
{Nova:1955}), or they are supposed to be trivial, but should be explicitly
writen down
here: ({\em i}) a region $S$ and an interval $\Delta t$ must be
prescribed at every point $(t,x^a)$ $\in {\cal E}$ in order
to define an averaged
field
$\langle p_\beta ^\alpha (t,x^a)\rangle _{{\cal E}}$, ({\em ii})
all the regions $S$ and time intervals $\Delta t$ are typical in some
defined sense - they are usually required to be of the same shape and volume,
$V_S={\rm const}$ and $
\Delta t={\rm const}$, and related to each other by
shifting along the Cartesian coordinate lines.
These properties are very easily arranged by
Lie-dragging of a region $S$ and an interval $\Delta t$, chosen around a point 
$(t,x^a)$ $\in {\cal E}$, along the congruences of the Cartesian coordinate
lines\footnote{
A definition of Lie-dragging, or dragging, of a region of a manifold along a
vector field (congruence), or in other words, mapping of a region into
another along a vector field (congruence) can be found in any standard
textbook on differential geometry (see, for example, \cite{Koba-Nomi:1963},
\cite{Schu:1980}). Throughout the paper Lie-draggings of regions are supposed 
to be $C^\infty$-diffeomorphisms.} to get a ``covering'' of the
manifold (or its connected part) with a region of the same shape and volume
and a time interval of the same length around each point of ${\cal E}$. 
Consequences of the
properties ({\em i}) and ({\em ii}) are, first of all,
the 
{\em uniqueness} of the definition of the average field $\langle p_\beta
^\alpha (t,x^a)\rangle _{{\cal E}}$, and, secondly, 
the commutation formulae between the averaging and the partial derivatives 
\begin{equation}
\label{comm:ED}\frac \partial {\partial t}\langle p_\beta ^\alpha
(t,x^a)\rangle _{{\cal E}}=\langle \frac \partial {\partial t}p_\beta
^\alpha (t,x^a)\rangle _{{\cal E}}\ ,\quad \frac \partial {\partial
x^a}\langle p_\beta ^\alpha (t,x^a)\rangle _{{\cal E}}=\langle \frac
\partial {\partial x^a}p_\beta ^\alpha (t,x^a)\rangle _{{\cal E}}\ .
\end{equation}
The essence of this approach lies in the exploitation of the calculational
advantages of the Cartesian coordinates on a flat manifold. The Cartesian
coordinates play a central role both in defining the
averages (\ref{defaver:ED}) and in obtaining properties, such as (\ref
{comm:ED}), which will allow the  
averaging out of the partial (ordinary) differential equations under interest.
It should be stressed that by requiring ({\em i}) and ({\em 
ii}) the averages $\langle p_\beta ^\alpha (t,x^a)\rangle _{
{\cal E}}$ become local functions
of $(t,x^a)$, i.e. 
\begin{equation}
\label{comm2:ED}\left( \frac \partial {\partial x^a}\frac \partial {\partial
t}-\frac \partial {\partial t}\frac \partial {\partial x^a}\right) \langle
p_\beta ^\alpha (t,x^a)\rangle _{{\cal E}}=0 \,, 
\end{equation}
and, therefore, they are provided with proper analytical properties.
Furthermore, the functional dependence of the average (\ref{defaver:ED})
on the averaging region and interval, 
$S$ and $\Delta t$, becomes simply a parametric dependence on
the value of the volume $V_S$ and the length $\Delta t${\rm . }One can,
therefore,
apply the standard differential and integral calculus to deal with the
averages and averaged equations. The set of averaged components $\langle
p_\beta ^\alpha (t,x^a)\rangle _{{\cal E}}$ given by (\ref{defaver:ED}) in
Cartesian coordinates are the components of a Lorentz tensor within the 
class of coordinate transformations 
$\widetilde{x}^\alpha =\Lambda _\beta ^\alpha x^\beta
+a^\beta $ with a constant shift $a^\beta $ and a constant Lorentz
transformation matrix $\Lambda _\beta ^\alpha $.
Now, if it is necessary
to consider an average field in another coordinate system, the
components of $\langle p_\beta ^\alpha (t,x^a)\rangle _{{\cal E}}$ in that
new system are found by applying the standard tensorial transformation
law.

There is still another property of the averages $\langle p_\beta ^\alpha
(t,x^a)\rangle _{{\cal E}}$ widely used in macroscopic electrodynamics 
(but again, often only tacitly assumed); the
idempotency of the averages 
\begin{equation}
\label{idem:ED}\langle \langle p_\beta ^\alpha (t,x^a)\rangle _{{\cal E}
}\rangle _{{\cal E}}=\langle p_\beta ^\alpha (t,x^a)\rangle _{{\cal E}}\ .
\end{equation}
In order to prove this property we must calculate 
the average value $\langle \langle p_\beta ^\alpha (t,x^a)\rangle _{{\cal E}
}\rangle _{{\cal E}}$ of the average $\langle p_\beta ^\alpha (t,x^a)\rangle
_{{\cal E}}$
\begin{equation}
\label{avav}\langle \langle p_\beta ^\alpha (t,x^a)\rangle _{{\cal E}
}\rangle _{{\cal E}}=\frac 1{\Delta tV_S}\int_{\Delta t}\int_S\left( \frac
1{\Delta t^{\prime }V_{S^{\prime }}}\int_{\Delta t^{\prime }}\int_{S^{\prime
}}p_\beta ^\alpha (t+t^{\prime }+t^{\prime \prime },x^a+x^{\prime
a}+x^{\prime \prime a})dt^{\prime \prime }d^3x^{\prime \prime }\right)
dt^{\prime }d^3x^{\prime }\ .
\end{equation}
Now, the expression (\ref{avav}) leads to (\ref{idem:ED}) under either of two
additional assumptions: ({\em iii}$^{\prime }$) the averaging region $\Delta
t^{\prime }\times S^{\prime }$ is the same region $\Delta t\times S$ with
the supporting point at $(t^{\prime },x^{\prime a})\in \Delta t\times S$ and
the average value $\langle p_\beta ^\alpha (t,x^a)\rangle _{{\cal E}}$ does
not depend on the choice of a supporting point within a chosen region; ({\em 
iii}$^{\prime \prime }$) the averaging region $\Delta t^{\prime }\times
S^{\prime }$ is a different neighbouring region
and the average value $\langle p_\beta ^\alpha (t,x^a)\rangle _{{\cal E}}
$ remains the same if evaluated over any neighbouring
regions containing a neighbourhood of $x$. Although they seem different the
above assumptions are nevertheless essentially equivalent. Indeed, ({\em iii}
$^{\prime }$) puts emphasis on the independence of the average value with
respect to the choice of a supporting point from a set of all possible
points in a {\em fixed averaging region}, while ({\em iii}$^{\prime
\prime }$) puts emphasis on the independence of the average value with
respect to the choice of an averaging region from a set of all possible
regions defined by a neighbourhood of a {\em fixed supporting point}.
Remembering ({\em ii}) that the averaging regions are typical, 
change of the supporting point can be considered as change of the
averaging region, and vice versa. Both assumptions thereby encode the same
fundamental property of the averages which is inherently related to the
philosophy of averaging itself - an averaging region is considered as a
point from the macroscopic point of view, and the change of either a
supporting point or an averaging region does not affect the corresponding
average value. The reason for this is that the differences 
caused by the above variations
are negligibly small from a macroscopic point of view\footnote{
It should be noted here that in classical hydrodynamics, as distinct from
macroscopic electrodynamics, discussion on the definition and properties
of the averages continues for more than one hundred years. A definition of an
average (either over space, time, ensemble, or a combination of such) in
hydrodynamics and its properties are vital elements of the theory itself for
it is clearly understood that the form of the equations  depends
on the definition and properties of the average. The definition (\ref
{defaver:ED}) under conditions ({\em i}) and ({\em ii}) with the properties 
(\ref{comm:ED}), (\ref{comm2:ED}) and (\ref{idem:ED}), which are part of
the Reynolds conditions in hydrodynamics, is known to result in the
fundamental equations of hydrodynamics describing the dynamics of 
turbulence. If one of
the Reynolds conditions is absent one must get different equations.
For a discussion on averages and their properties in hydrodynamics, see,
for example, \cite{Moni-Yagl:1965} and references therein.}.

Thus, the well-known procedure of space-time averaging in classical
electrodynamics presupposes the above specific conditions ({\em i}), 
({\em ii}) and ({\em iii}) to ensure reasonable analytical and 
tensorial properties
of the averages (\ref{defaver:ED}). Clearly, the whole procedure 
relies crucially  
on the existence and properties of the exceptional
coordinate system in a flat space-time, namely, Cartesian coordinates. These
conditions, which seem more or less trivial at first sight, 
require further analysis to make clear their geometrical meaning and
invariant content. This will allow a reasonable formulation of a 
space-time
averaging procedure on general (not necessarily flat) manifolds.

\section{Definition of space-time averages in macroscopic gravity}

Let us remind the definition of the space-time averages adopted in
macroscopic gravity \cite{Zala:1992},\cite{Zala:1993}. This 
procedure is a generalization of the space-time averaging
procedure adopted in electrodynamics (see Section 2) and it is also based on
the concept of Lie-dragging of averaging regions, which makes it valid for
any differentiable manifold.

\begin{definition}
\label{average} Chosen a compact region $\Sigma \subset {\cal M}$ in an 
$n$-dimensional differentiable metric manifold 
$({\cal M}$, $g_{\alpha \beta })$ with a volume $n$-form 
and a supporting point $x\in \Sigma $ to which the average value will be
prescribed, the average value of an object (tensor, geometric object,
etc.) $p_\beta ^\alpha (x),\,x\in {\cal M\ }$, over a region $\Sigma $ at
the supporting point $x\in \Sigma $ is defined as 
\begin{equation}
\label{defaver:MG}\overline{p}_\beta ^\alpha (x)=\frac 1{V_\Sigma
}\int_\Sigma {\bf p}_\beta ^\alpha (x,x^{\prime })\sqrt{-g^{\prime }}
d^nx^{\prime }\equiv \langle {\bf p}_\beta ^\alpha \rangle \ ,
\end{equation}
where $V_\Sigma $ is the volume of the region $\Sigma $, 
\begin{equation}
\label{volume}V_\Sigma =\int_\Sigma \sqrt{-g}d^nx\ .
\end{equation}
\end{definition}

Here the integration is carried out over all points 
$x^{\prime }\in \Sigma $, $g^{\prime }=\det (g_{\alpha \beta }(x^{\prime }))$
(since the primary interest is in space-time manifolds, it is
assumed that this 
determinant is negative, otherwise
the negative sign in $\sqrt{-g}$ must be dropped)
and the bold face object ${\bf p}_\beta
^\alpha (x,x^{\prime })$ in the integrand of (\ref{defaver:MG}) is a bilocal
extension of the object $p_\beta ^\alpha (x)$, 
\begin{equation}
\label{bilocext}{\bf p}_\beta ^\alpha (x,x^{\prime })={\cal A}_{\mu ^{\prime
}}^\alpha (x,x^{\prime })p_{\nu ^{\prime }}^{\mu ^{\prime }}(x^{\prime })
{\cal A}_\beta ^{\nu ^{\prime }}(x^{\prime },x)\ , 
\end{equation}
by means of bilocal averaging operators ${\cal A}_{\beta ^{\prime }}^\alpha
(x,x^{\prime })$ and ${\cal A}_\beta ^{\alpha ^{\prime }}(x^{\prime },x)$.
The averaging scheme is covariant and linear, $\langle a
{\bf p}_\beta ^\alpha +b{\bf q}_\beta ^\alpha
\rangle =a\langle {\bf p}_\beta ^\alpha \rangle +b\langle {\bf q}_\beta
^\alpha \rangle $,  $a,b\in {\rm R}$, 
by construction and the averaged object $
\overline{p}_\beta ^\alpha $ keeps the same tensorial character as $
p_\beta ^\alpha $.

Let us suppose that the bilocal functions ${\cal A}_{\beta ^{\prime
}}^\alpha (x,x^{\prime })$ and ${\cal A}_\beta ^{\alpha ^{\prime
}}(x^{\prime },x)$ are defined locally on an open subset
${\cal U}\subset {\cal M}$, $
x,x^{\prime }\in {\cal U}$. In the following sections
their algebraic and differential
properties are formulated and analyzed to show that averaging
operators with such properties do exist and also to find out the corresponding
properties of the averages (\ref{defaver:MG}).

\section{Algebra of the averaging operator}

The following algebraic properties, which are a formalization of the
properties of the space-time averages in macroscopic electrodynamics using the
language of bilocal operators, are required to hold \cite{Zala:1996},\cite
{Zala:1992},\cite{Zala:1993}.

\begin{property}
\label{coincidence1} The coincidence limit of ${\cal A}_{\beta ^{\prime
}}^\alpha $ is 
\begin{equation}
\label{coincidence}\lim _{x^{\prime }\rightarrow x}{\cal A}_{\beta ^{\prime
}}^\alpha (x,x^{\prime })=\delta _\beta ^\alpha \ .
\end{equation}
\end{property}

\begin{property}
\label{idempotency1} The operator ${\cal A}_{\beta ^{\prime }}^\alpha $ is
idempotent 
\begin{equation}
\label{idempotency}{\cal A}_{\beta ^{\prime }}^\alpha (x,x^{\prime }){\cal A}
_{\gamma ^{\prime \prime }}^{\beta ^{\prime }}(x^{\prime },x^{\prime \prime
})={\cal A}_{\gamma ^{\prime \prime }}^\alpha (x,x^{\prime \prime })\ .
\end{equation}
\end{property}

These two properties imply that ${\cal A}_\beta ^{\alpha
^{\prime }}(x^{\prime },x)$ is the inverse operator of ${\cal A}_{\beta ^{\prime
}}^\alpha (x,x^{\prime })$, ${\cal A}_{\beta ^{\prime }}^\alpha {\cal A}
_\gamma ^{\beta ^{\prime }}=\delta _\gamma ^\alpha $ and ${\cal A}_{\beta
^{\prime }}^\alpha {\cal A}_\alpha ^{\gamma ^{\prime }}=\delta _{\beta
^{\prime }}^{\gamma ^{\prime }}$, and that the average tensor $\overline{p}_\beta
^\alpha (x)$ takes the same value as the original tensor $p_\beta ^\alpha
(x) $, $\overline{p}_\beta ^\alpha (x)=p_\beta ^\alpha (x)$,
when the integrating region $\Sigma $ is chosen infinitesimally small, or
tends to zero. This implies an additional algebraic property of the
averages (\ref{defaver:MG}), namely, that the averaging procedure commutes with
the operation of index contraction.

The idempotency (\ref{idempotency}) of the averaging operator ${\cal A}
_{\beta ^{\prime }}^\alpha $ is designed to provide the idempotency of the
averages (\ref{defaver:MG}) in macroscopic gravity 
\begin{equation}
\label{aver=aver}\overline{\overline{p}}_\beta ^\alpha (x)=\overline{p}
_\beta ^\alpha (x)\ .
\end{equation}
Indeed, let us consider for simplicity a vector $v^\alpha (x)$ and calculate
the twice averaged value $\overline{\overline{v}}^\alpha (x)$ around the same
point $x\in {\cal U}$ by using the definition (\ref{defaver:MG}) 
\begin{equation}
\label{avravr}\overline{\overline{v}}^\alpha (x)=\frac 1{V_{\Sigma
_x}}\int_{\Sigma _x}\left( \frac 1{V_{\Sigma _{x^{\prime }}}}\int_{\Sigma
_{x^{\prime }}}{\cal A}_{\beta ^{\prime }}^\alpha (x,x^{\prime }){\cal A}
_{\gamma ^{\prime \prime }}^{\beta ^{\prime }}(x^{\prime },x^{\prime \prime
})v^{\gamma ^{\prime \prime }}(x^{\prime \prime })\sqrt{-g^{\prime \prime }}
d^nx^{\prime \prime }\right) \sqrt{-g^{\prime }}d^nx^{\prime }
\end{equation}
where $\Sigma _{x^{\prime }}$ is an averaging region around the point 
$x^{\prime }\in {\cal U}$. By applying the idempotency condition (\ref
{idempotency}) the expression (\ref{avravr}) takes the form 
\begin{equation}
\label{vbarbar}\overline{\overline{v}}^\alpha (x)=\frac 1{V_{\Sigma
_x}}\int_{\Sigma _x}\left( \frac 1{V_{\Sigma _{x^{\prime }}}}\int_{\Sigma
_{x^{\prime }}}{\cal A}_{\gamma ^{\prime \prime }}^\alpha (x,x^{\prime
\prime })v^{\gamma ^{\prime \prime }}(x^{\prime \prime })\sqrt{-g^{\prime
\prime }}d^nx^{\prime \prime }\right) \sqrt{-g^{\prime }}d^nx^{\prime }\ .
\end{equation}
Now, if the term in parentheses were independent of $x^{\prime }$, we could
take it outside the integral with respect to the variable $x^{\prime }$ to
get the resulting property (\ref{aver=aver}). This term in (\ref{vbarbar}),
however, depends explicitly on $x^{\prime }$ in the integration region $
\Sigma _{x^{\prime }}$. Similarly as in macroscopic electrodynamics, 
the idempotency (\ref{aver=aver}) of 
the averages (\ref{defaver:MG}) follows under either of
two additional assumptions \cite{Zala:unpub}: ({\em iii}$^{\prime }$) the
averaging region $\Sigma _{x^{\prime }}$ is the same region $\Sigma _x$ with
the supporting point at $x^{\prime }\in \Sigma _x$ and the average value $
\overline{p}_\beta ^\alpha (x)$ does not depend on the choice of a
supporting point within a chosen region; ({\em iii}$^{\prime \prime }$) the
averaging region $\Sigma _{x^{\prime }}$ is a different neighbouring region
and the average value $\overline{p}_\beta ^\alpha (x)$ remains the same
if evaluated over any neighbouring region containing a
neighbourhood of $x$. As it was emphasized in Section 2 this is a
fundamental property of any physically reasonable classical averaging
procedure and such stability of averages comes from the basic principles of
averaging. 

A microscopic field to be averaged is supposed to have two
essentially different variation scales\footnote{
The case of two scales is discussed here for the sake of simplicity. Of
course, very often there is a hierarchy of scales, in which case the
arguments are applied for each couple of scales satisfying (\ref{scales}) 
to be micro- and macroscopic
ones, respectively.}, $\lambda $ and $L$, satisfying 
\begin{equation}
\label{scales}
\lambda <<L\ , 
\end{equation}
and an averaging region must be taken of an intermediate size $d$ such as 
\begin{equation}
\label{region}
\lambda <<d<<L
\end{equation}
so that the averaging effectively smooths out all the variations of the
microscopic field of the scale $\lambda $. It is implicitly assumed in every
averaging scheme (for example, in the averaging schemes applied in classical
physics in Minkowski space-time, such as in hydrodynamics and
electrodynamics - see Section 2) that the result of the averaging is
insensitive to the choice of the supporting point within a fixed averaging
region and it is independent of the choice of integration (averaging) region
itself provided the scale $d$ satisfies the condition (\ref{region}).
This means that the microscopic averaging region is considered as a single
``point'' for the macroscopic averaged field. Such regions have been called
``physically infinitesimally small'' by Lorentz \cite{Lore:1916}.

The trouble is, however, to perform a mathematically rigorous proof of this
fact. This problem already appears in the usual averaging procedures in flat 
space-times and,
to the best of our knowledge, remains unsolved (see, for example, a discussion
in \cite{Moni-Yagl:1965} for the case of hydrodynamics). Furthermore,
it is not clear
whether the property (\ref{aver=aver}) has to be interpreted as an exact
one, or it is rather an approximate property. A satisfactory formal
analysis of the idempotency property of the volume averages is still
lacking, though it is extensively used\footnote{
It should be pointed out here that within the ensemble average procedure
the
idempotency follows without problem. However, such averages have
their own problems. In particular, the whole body of 
problems related with idempotency is replaced by the necessity to
prove the ergodicity hypothesis which states that ensemble and time (or
space) averages are equivalent. Both ensemble and volume averagings 
have their own advantages and areas of applicability in describing
physical phenomena. It is important to realize in this connection 
that in all macroscopic
settings a volume averaging (over space, time, or space-time) is an
unavoidable element \cite{deGr:1969}-\cite{Jack:1975}.}. This problem
certainly deserves further examination \cite{Mars-Zala:progr}.

\section{Structure of the averaging operator}

In the previous section we discussed two algebraic Properties 
\ref{coincidence1} and \ref{idempotency1} of the bilocal operator ${\cal A}
_{\beta ^{\prime }}^\alpha (x,x^{\prime })$ and the physical motivations
which make these two conditions plausible. In this section the consequences
of requiring (\ref{coincidence}) and (\ref{idempotency}) in the structure of 
the bilocal operator 
are analyzed. The following main
theorem holds. 

\begin{theorem}
\label{factoriz} The bilocal operator ${\cal A}_\beta ^{\alpha ^{\prime
}}(x^{\prime },x)$ is idempotent if and only if it is factorized. That is to
say, for a bilocal operator ${\cal A}_\beta ^{\alpha ^{\prime }}(x^{\prime
},x)$ Properties \ref{coincidence1} and \ref{idempotency1} are
equivalent to 
\begin{equation}
\label{factor}{\cal A}_\beta ^{\alpha ^{\prime }}(x^{\prime },x)=F_i^{\alpha
^{\prime }}(x^{\prime })F{^{-1}}_\beta ^i(x)
\end{equation}
where $F_i^\beta (x)$ is a set of $n$ linear 
independent vector fields and $F{^{-1}}_\beta ^i(x)$ is the 
associated dual 1-form basis, $i=1,...,n$, and
summation over $i$ is carried out.
\end{theorem}

{\it Proof.} The inverse implication is checked trivially,
so let us consider the direct implication only. Indeed, the partial
derivative of  the idempotency relation (\ref{idempotency1}), ${\cal A}
_{\alpha ^{\prime }}^{\gamma ^{\prime \prime }}(x^{\prime \prime },x^{\prime
}){\cal A}_\beta ^{\alpha ^{\prime }}(x^{\prime },x)={\cal A}_\beta ^{\gamma
^{\prime \prime }}(x^{\prime \prime },x)$,${\cal \ }$with respect to $
x^{\sigma ^{\prime }}$ results in the following expression 
$$
{\cal A}_{\alpha ^{\prime} , \sigma ^{\prime }}^{\gamma ^{\prime
\prime }}(x^{\prime \prime },x^{\prime }){\cal A}_\beta ^{\alpha ^{\prime
}}(x^{\prime },x)+{\cal A}_{\alpha ^{\prime }}^{\gamma ^{\prime \prime
}}(x^{\prime \prime },x^{\prime }){\cal A}
_{\beta , \sigma ^{\prime }}^{\alpha ^{\prime }}(x^{\prime },x)=0\ . 
$$
Contracting this expression with ${\cal A}_{\gamma ^{\prime \prime
}}^{\varepsilon ^{\prime }}(x^{\prime },x^{\prime \prime })$ and ${\cal A}
_{\mu ^{\prime }}^\beta (x,x^{\prime })$ and using both (\ref{coincidence})
and (\ref{idempotency}) one gets 
\begin{equation}
{\cal A}_{\gamma ^{\prime \prime }}^{\epsilon ^{\prime }}(x^{\prime
},x^{\prime \prime }){\cal A}_{\mu ^{\prime} , \sigma ^{\prime }}^{\gamma 
^{\prime \prime }}(x^{\prime \prime },x^{\prime })+
{\cal A}_{\beta , \sigma ^{\prime }}^{\epsilon ^{\prime }}(x^{\prime },x)
{\cal A}_{\mu ^{\prime }}^\beta (x,x^{\prime })=0\ .
\end{equation}
Now, the first summand does not depend on $x$ and the second one does not
depend on $x^{\prime \prime }$, which means that 
\begin{equation}
\label{defB}{\cal A}_{\beta , \sigma ^{\prime }}^{\epsilon ^{\prime
}}(x^{\prime },x){\cal A}_{\mu ^{\prime }}^\beta (x,x^{\prime })\equiv
B_{\mu ^{\prime }\sigma ^{\prime }}^{\epsilon ^{\prime }}(x^{\prime }),
\end{equation}
where $B_{\mu ^{\prime }\sigma ^{\prime }}^{\alpha ^{\prime }}$ are
arbitrary functions depending only on the variable $x^{\prime }$. The
relation (\ref{defB}) immediately implies the following identity for the
functions $B_{\mu ^{\prime }\sigma ^{\prime }}^{\alpha ^{\prime }}(x^{\prime
})$: 
\begin{equation}
\label{intrB}B_{\mu ^{\prime }\left[ \sigma ^{\prime },\delta ^{\prime
}\right] }^{\epsilon ^{\prime }}+B_{\rho ^{\prime }[\sigma ^{\prime
}}^{\epsilon ^{\prime }}B_{\underline{\mu ^{\prime }}\delta ^{\prime
}]}^{\rho ^{\prime }}=0\ ,
\end{equation}
where the square brackets denote, as
usual, antisymmetrization and the underlined indices are not affected by
antisymmetrization. Considering now (\ref{defB}) as a set of linear partial
differential equations for the unknowns ${\cal A}_\beta ^{\alpha ^{\prime
}}(x^{\prime },x)$ ($x$ being a parameter)
\begin{equation}
\label{equA}{\cal A}_{\beta , \sigma ^{\prime }}^{\epsilon ^{\prime
}}(x^{\prime },x)=B_{\mu ^{\prime }\sigma ^{\prime }}^{\epsilon ^{\prime
}}(x^{\prime }){\cal A}_\beta ^{\mu ^{\prime }}(x^{\prime },x)\ ,
\end{equation}
one observes that they are always integrable because $B_{\mu
^{\prime }\sigma ^{\prime }}^{\alpha ^{\prime }}(x^{\prime })$ satisfy the
identity (\ref{intrB}). The general solution of the system of
equations (\ref{equA})
is of the form
$$
{\cal A}_\beta ^{\epsilon ^{\prime }}(x^{\prime },x)=F_i^{\epsilon ^{\prime
}}(x^{\prime })H_\beta ^i(x) 
$$
where $F_i^{\alpha ^{\prime }}(x^{\prime })$, $i=1,...,n$, are $n^2$
independent solutions of (\ref{equA}) and $H_\beta ^i(x)$ are $
n^2$ arbitrary constants of integration (thereby depending on $x$). It only
remains to impose the coincidence limit property (\ref{coincidence}) to set $
H_\beta ^i(x)=F{^{-1}}_\beta ^i(x)$ and the proof of the theorem is
completed. $\quad \Box $

\section{Differential properties of bilocal operators}

In this section the differential properties
of the bilocal operators are summarized. This overview is essential 
for understanding further results concerning the space-time averaging
scheme.

In order to obtain the averaged fields of the geometric objects on ${\cal M}$
one needs to assign an averaging region $\Sigma _x$ to each point $x$ of $
{\cal U}\subset {\cal M}$, where the averaging integral (\ref{defaver:MG}) is
to be evaluated. Furthermore, to calculate directional, partial and
covariant derivatives of the averaged fields a law of correspondence
between neighbouring averaging regions must be defined \cite{Zala:1992},\cite
{Zala:1993},\cite{AZS:1990b}. A natural  way to define such a
correspondence is to relate averaging regions by Lie-dragging, or mapping
of a region into another along a vector field (see Section 2), by means of
another bilocal operator ${\cal W}_\beta ^{\alpha ^{\prime }}(x^{\prime },x)$
which is also assumed to satisfy the coincidence limit property (\ref
{coincidence}).

To derive the commutation formulae between
the averaging and the derivation, one
should define first the directional derivative of an average field $\overline{p
}_\beta ^\alpha (x)$ along a vector field $\vec \xi =d/d\lambda $, 
\begin{equation}
\label{derivative}\frac d{d\lambda }\overline{p}_\beta ^\alpha (x)=\lim
_{\Delta \lambda \rightarrow 0}\frac 1{\Delta \lambda }[\overline{p}_\beta
^\alpha (x+\Delta x)-\overline{p}_\beta ^\alpha (x)] 
\end{equation}
where $\Delta x^\alpha =\xi ^\alpha (x)\Delta \lambda $. Let us define now
the shift field for every point $x^{\prime }\in \Sigma _x$
as \cite{Zala:1992},\cite{Zala:1993},\cite{AZS:1990b} 
\begin{equation}
\label{shift}S^{\alpha ^{\prime }}(x^{\prime },x)={\cal W}_\beta ^{\alpha
^{\prime }}(x^{\prime },x)\xi ^\beta (x)\ .
\end{equation}
Now,  the averaging region $\Sigma _{x+\Delta x}$ associated with
the point $x+\Delta x$ is obtained by Lie-dragging the averaging region
$\Sigma _x$ a parametric length 
$\Delta \lambda $ (the same for all $x^{\prime }\in
\Sigma _x$) along the integral lines of the field $S^{\alpha ^{\prime }}$.

The coordination  bivector ${\cal W}_\beta ^{\alpha ^{\prime }}$
allows the construction of the shift vector $S^{\alpha ^{\prime }}$ for
any averaging region and any vector $\xi ^\alpha $ (the
shift vector at the supporting point for that region). By choosing $n$ such
linearly independent vector fields $\xi _i^\alpha $ and shifting averaging
regions along them one can build a covering of the manifold with one 
averaging region associated to every
$x\in {\cal U}\subset {\cal M}$.
This procedure is a formalization of the condition ({\em i}) of
Section 2.

As a consequence of the definitions (\ref{defaver:MG}), (\ref{derivative})
and (\ref{shift}) one can obtain the following formula for the
commutation of partial
differentiation and averaging \cite{Zala:1992},\cite{Zala:1993}: 
\begin{equation}
\label{derivative2}\overline{p}_{\beta ,\lambda }^\alpha =\langle {\cal A}
_{\mu ^{\prime }}^\alpha p_{\nu ^{\prime },\epsilon ^{\prime }}^{\mu
^{\prime }}{\cal A}_\beta ^{\nu ^{\prime }}{\cal W}_\lambda ^{\epsilon
^{\prime }}\rangle +\langle {\bf p}_\beta ^\alpha {\cal W}_{\lambda
:\epsilon ^{\prime }}^{\epsilon ^{\prime }}\rangle -\overline{p}_\beta
^\alpha \langle {\cal W}_{\lambda :\epsilon ^{\prime }}^{\epsilon ^{\prime
}}\rangle -\langle {\cal S}_{\sigma \lambda }^\alpha {\bf p}_\beta ^\sigma
\rangle +\langle {\bf p}_\sigma ^\alpha {\cal S}_{\beta \lambda }^\sigma
\rangle \ . 
\end{equation}
Here, ${\cal S}_{\sigma \lambda }^\alpha $ are the so-called structural
functions, ${\cal S}_{\beta \gamma }^\alpha ={\cal A}_{\epsilon ^{\prime
}}^\alpha ({\cal A}_{\beta ,\gamma }^{\epsilon ^{\prime }}+{\cal A}_{\beta
,\sigma ^{\prime }}^{\epsilon ^{\prime }}{\cal W}_\gamma ^{\sigma ^{\prime
}})$ and ${\cal W}_{\lambda :\epsilon ^{\prime }}^{\epsilon ^{\prime }}$ is
the divergence of the coordination bivector, ${\cal W}_{\lambda :\epsilon
^{\prime }}^{\epsilon ^{\prime }}={\cal W}_{\lambda ,\epsilon ^{\prime
}}^{\epsilon ^{\prime }}+(\ln \sqrt{-g})_{,\epsilon ^{\prime }}{\cal W}
_\lambda ^{\epsilon ^{\prime }}$. In an affine connection space with 
connection coefficients $\Gamma^\alpha_{\beta \gamma}$ (a Riemannian space 
is considered in \cite{AZS:1990b}) the same
formula (\ref{derivative2}) remains valid for covariant
derivatives, where partial differentiation is replaced by the covariant
one and the divergence ${\cal W}_{\lambda
:\epsilon ^{\prime }}^{\epsilon ^{\prime }}$ is replaced by ${\cal W}
_{\lambda ;\epsilon ^{\prime }}^{\epsilon ^{\prime }}={\cal W}_{\lambda
,\epsilon ^{\prime }}^{\epsilon ^{\prime }}+\Gamma _{\epsilon ^{\prime }}
{\cal W}_\lambda ^{\epsilon ^{\prime }}$ 
(with  $\Gamma _{\epsilon ^{\prime }}=
\Gamma ^{\alpha^{\prime}}_{\epsilon^{\prime } \alpha^{\prime }}$). 
The commutation
formula (\ref{derivative2}) 
has a very transparent meaning: the first term in the right-hand side
is the average value of the
derivative of $p_\beta ^\alpha $ weighted by ${\cal W}_\beta ^{\alpha ^{\prime
}}$, the last two terms are due to the non-triviality of the
averaging operator ${\cal A}_\beta ^{\alpha ^{\prime }}$, while the second
and third terms describe
the effect of a non-trivial averaging measure
in (\ref{defaver:MG}) and the variation in the value of the volume for
different
regions, respectively. Indeed, the change in the volume $V_\Sigma $ of an (averaging)
region $\Sigma $ $\subset {\cal M}$ Lie-dragged along a vector field $\xi $
is given by (see, for example, \cite{AZS:1990a},\cite{Schu:1980}) 
\begin{equation}
\label{volchange}\frac d{d\lambda }V_\Sigma =\int_\Sigma {\rm div}\,\xi
d\Omega \ , 
\end{equation}
which can be written in terms of partial derivatives after using (\ref{shift})
as \cite{Zala:unpub}
\begin{equation}
\label{shift:volume}V_{\Sigma ,\beta }=\langle {\cal W}_{\beta :\alpha
^{\prime }}^{\alpha ^{\prime }}\rangle V_\Sigma \,\ . 
\end{equation}
The expression (\ref{derivative2}) is the most general version of formulae 
(\ref{comm:ED}) for arbitrary averaging and coordination operators.

A fundamental problem one has to face now is whether it is possible or not
to define a unique covering of the manifold ${\cal M}$, like
in Minkowski space-time (see Section 2). In other words, given a microscopic
tensor field $p_\beta ^\alpha (x)$, the goal is to determine a uniquely
defined averaged tensor field $\overline{p}_\beta ^\alpha (x)$ with
reasonable analytical properties in its dependence on the supporting point.
Furthermore, we must deal with the problem that
the averages (\ref{defaver:MG}) depend functionally on the
averaging regions $
\overline{p}_\beta ^\alpha (x)=\overline{p}_\beta ^\alpha (x)[\Sigma ]$ 
and applying  the commutation formula (\ref{derivative2}) to find
averaged equations would, in general,
bring volume dependent terms into them, which is undesirable.
In the case of Minkowski space-time the uniqueness of the
averaged field with the proper analytical behaviour (\ref{comm2:ED}) and the
parametric dependence on the averaging region volume are ensured by the
condition ({\em ii}) (see Section 2) which sets the averaging regions to be
typical, of the same shape and volume, 
$V_S={\rm const}$ and $\Delta t={\rm const}$, and shifted along the
Cartesian coordinate lines. To arrange  similar properties for the generalized
averages (\ref{defaver:MG})
one must look for some specific conditions on the coordination operator $
{\cal W}_\beta ^{\alpha ^{\prime }}$. The following remarkable
theorem holds \cite{Zala:1992},\cite{Zala:1993} (its version for Riemannian
manifolds has been given in \cite{AZS:1990b}).

\begin{theorem}
\label{uniqueness} In the averaging region coordination by the bivector
${\cal W}_\beta
^{\alpha ^{\prime }}$ on an arbitrary differentiable manifold, 
it is necessary and sufficient to require 
\begin{equation}
\label{diffW:1}{\cal W}_{[\beta ,\gamma ]}^{\alpha ^{\prime }}+{\cal W}
_{[\beta ,\delta ^{\prime }}^{\alpha ^{\prime }}{\cal A}_{\gamma ]}^{\delta
^{\prime }}=0\,\ ,
\end{equation}
for the average tensor field $\overline{p}_\beta ^\alpha (x)$ to be a
single valued local function of the supporting point $x$ on ${\cal U}\subset 
{\cal M}$, 
\begin{equation}
\label{averlocal}\overline{p}_{\beta ,[\mu \nu ]}^\alpha =0\ .
\end{equation}
\end{theorem}

{\it Proof.} The sufficiency of (\ref{diffW:1}) is easily proved
by calculating the antisymmetrized 
second partial derivative of (\ref{derivative2}) which, after
using  (\ref{diffW:1}), gives (\ref{averlocal}). The necessity follows from
the analysis of the condition (\ref
{averlocal}) after taking into account that the averaging regions and the
tensor
field $p_\beta ^\alpha $ are arbitrary. $\quad \Box $

Geometrically, (\ref{diffW:1}) means that, given an averaging region $\Sigma 
$, the region $\widetilde{\Sigma }$ obtained by transporting 
$\Sigma $ along an infinitesimal parallelogram constructed from two
commuting vector fields $\xi $ and $\zeta $ according
to the law (\ref{shift}) coincides with the original region, $
\widetilde{\Sigma }=\Sigma $. This is a highly non-trivial property
which allows to construct a covering of the manifold with an averaging region
attached to every point in the manifold, thus generalizing the
corresponding part of the condition ({\em ii}) in Section 2.
In the formalism of bilocal exterior calculus the condition (\ref{diffW:1})
reads that the operator ${\cal W}_\beta ^{\alpha ^{\prime }}$ is
biholonomic, which means that the bilocal coordinate 1-form basis ${\cal W}
^{\alpha ^{\prime }}={\cal W}_\beta ^{\alpha ^{\prime }}{\rm d}x^\beta $ has
vanishing biholonomicity coefficients in the bilocal Maurer-Cartan equations 
\cite{Zala:1993}.

It should be noted here that formula (\ref{averlocal}) is analogous to
formula (\ref{comm2:ED}) in both form and meaning.

Another condition on the coordination bivector ${\cal W}_\beta ^{\alpha
^{\prime }}$ is the requirement  that the Lie-dragging of a region is 
a volume-preserving diffeomorphism \cite{Zala:1992},\cite{Zala:1993}
(\cite{AZS:1990b} for Riemannian manifolds) 
\begin{equation}
\label{diffW:2}{\cal W}_{\beta :\alpha ^{\prime }}^{\alpha ^{\prime }}=0\ ,
\end{equation}
which means that the averaging regions do not change the value of the
volume when shifted (coordinated) along a chosen vector field $\xi $ 
according to
(\ref{shift}). This generalizes the corresponding part of the condition 
({\em ii}) in Section 2.

Thus, the condition (\ref{diffW:1}) states that the average tensor field is
a single
valued local function of the supporting point $x$ (\ref{averlocal}). 
Adding the condition (\ref{diffW:2}), the average tensor field does not
depend explicitly on the value of the region volume $V=V_\Sigma $, and $V$
itself is a free parameter of the theory. Given a microscopic tensor
field $p_\beta ^\alpha (x)$ on ${\cal M}$, the average tensor field 
$\overline{p}_\beta ^\alpha (x)$ is
therefore uniquely defined on ${\cal U}\subset {\cal M}$ and can be handled
within the framework of standard differential and integral calculus. 
Requiring additionally that the two bivectors $
{\cal A}_\beta ^{\alpha ^{\prime }}$ and ${\cal W}_\beta ^{\alpha
^{\prime }}$ coincide
 \begin{equation}
\label{AequalW}{\cal A}_\beta ^{\alpha ^{\prime }}={\cal W}_\beta ^{\alpha
^{\prime }}\,
\end{equation}
the first term in the commutation formula (\ref{derivative2}) becomes
exactly the average derivative. Using all conditions (\ref{diffW:1}), 
(\ref{diffW:2}) and (\ref{AequalW}), the commutation
formula acquires a remarkable simply form \cite{Zala:1992},\cite
{Zala:1993} 
\begin{equation}
\label{comm}\overline{p}_{\beta ,\gamma }^\alpha =\langle {\bf p}_{\beta
,\gamma }^\alpha +{\bf p}_{\beta ,\alpha ^{\prime }}^\alpha {\cal W}^{\alpha
^{\prime }}_\gamma \rangle \ .\,
\end{equation}
The corresponding analogues of this expression
for covariant differentiation are obtained
by replacing partial derivatives by covariant ones. To obtain the
expression for
directional derivatives we must contract this expression  with a vector $\xi
^\gamma $ and insert the vector field $S^{\alpha ^{\prime }}$ from 
(\ref{shift}) in the second term of the right-hand side of (\ref{comm}). 
Formula (\ref{comm}) generalizes formulae (\ref{comm:ED}) and 
it can be easily shown \cite{Zala:1993} to become exactly (\ref{comm:ED}) if $
{\cal W}_\beta ^{\alpha ^{\prime }}=\delta _\beta ^\alpha $ and the volume $
n$-form $\varepsilon $ is standard, that is 
$(\ln \sqrt{-g})_{,\epsilon ^{\prime }}=0
$ and $\varepsilon = {\rm d}x^1 \wedge ... \wedge {\rm d}x^n$ 
\cite{Koba-Nomi:1963}, \cite{Schu:1980}. In Section 8 one can 
find more details on this particular case, see formulae 
(\ref{defaver:MG:prop}) and (\ref{commu:2}). 

\section{Existence theorems}

The differential conditions (\ref{diffW:1}) and (\ref{diffW:2}) together
with the algebraic conditions (\ref{coincidence}), (\ref{idempotency}) and 
(\ref{AequalW}) are to be considered as a set of partial differential and
algebraic equations for the unknown functions ${\cal W}_\beta ^{\alpha ^{\prime
}}$. Provided a solution for the system is found, the existence of such
operators ${\cal A}_\beta ^{\alpha ^{\prime }}$ and ${\cal W}_\beta ^{\alpha
^{\prime }}$ (and therefore of the averages with the above
described properties) is proved.
Theorem \ref{factoriz} has revealed the structure of the operator ${\cal W}
_\beta ^{\alpha ^{\prime }}$ obeying the algebraic properties (\ref
{coincidence}), (\ref{idempotency}). The following theorem gives the
general solution of (\ref{diffW:1}) 
(in \cite{Zala:1992},\cite{Zala:1993} and \cite {AZS:1990b} for Riemannian 
manifolds, the same theorem has proved a solution of (\ref{diffW:1}) with
(\ref{AequalW}) sought in a factorized form (\ref{factor}), and now,
with Theorem \ref{factoriz} taken into account, 
it gives the general solution). 
\begin{theorem}
\label{soldiffW:1} In an arbitrary $n$-dimensional differentiable manifold
the general solution of the equations 
\begin{equation}
\label{diffWW:1}{\cal W}_{[\beta ,\gamma ]}^{\alpha ^{\prime }}
+{\cal W}_{[\beta ,\delta ^{\prime }}^{\alpha ^{\prime }}
{\cal W}_{\gamma ]}^{\delta ^{\prime }}=0\,\ ,
\end{equation}
for idempotent bilocals ${\cal W}_\beta ^{\alpha ^{\prime }}(x^{\prime},x)$ 
is given by 
\begin{equation}
\label{W}{\cal W}_\beta ^{\alpha ^{\prime }}(x^{\prime },x)=f_i^{\alpha
^{\prime }}(x^{\prime })f{^{-1}}_\beta ^i(x)
\end{equation}
where $f_i^\alpha (x)\partial _\alpha =\mbox{\boldmath$f_{i}$}$ is any
vector basis satisfying the commutation relations 
\begin{equation}
\label{commff}\left[ \mbox{\boldmath$f_{i}$},\mbox{\boldmath$f_{j}$}\right]
=C_{ij}^k\mbox{\boldmath$f_{k}$}
\end{equation}
with constant structure functions (anholonomicity coefficients) $C_{ij}^k$, 
\begin{equation}
\label{C:const}C_{ij}^k={\rm const}.
\end{equation}
\end{theorem}

{\it Proof. }Due to Theorem \ref{factoriz} the idempotent
bivector ${\cal W}_\beta ^{\alpha ^{\prime }}(x^{\prime },x)$ has the
factorized form (\ref{factor})
$$
{\cal W}_\beta ^{\alpha ^{\prime }}(x^{\prime },x)=F_i^{\alpha ^{\prime
}}(x^{\prime })F{^{-1}}_\beta ^i(x)\,\ . 
$$
Inserting this expression into (\ref{diffWW:1}) one finds 
\begin{equation}
\label{diffFF:1}F_i^{\alpha ^{\prime }}(x^{\prime })F{^{-1}}_{[\beta ,\gamma
]}^i(x)+F_{i,\delta ^{\prime }}^{\alpha ^{\prime }}(x^{\prime })F{^{-1}}
_{[\beta }^i(x)F_j^{\delta ^{\prime }}(x^{\prime })F{^{-1}}_{\gamma
]}^j(x)=0\,\ . 
\end{equation}
Now, using the expression relating derivatives of a vector basis and its dual
1-form basis 
\begin{equation}
F_{j,\mu ^{\prime }}^{\alpha ^{\prime }}=-F_l^{\alpha ^{\prime }}F{^{-1}}
_{\beta ^{\prime },\mu ^{\prime }}^lF{^{-1}}_j^{\beta ^{\prime }}\,\ 
\end{equation}
in equation (\ref{diffFF:1}), it becomes 
\begin{equation}
\label{anhol:const}2F_l^\rho F_k^\sigma F{^{-1}}_{\left[ \rho ,\sigma
\right] }^i(x)=2F_l^{\rho ^{\prime }}F_k^{\sigma ^{\prime }}F{^{-1}}_{\left[
\rho ^{\prime },\sigma ^{\prime }\right] }^i(x^{\prime })\equiv C_{lk}^i\,\ , 
\end{equation}
which is exactly the expression for the anholonomicity coefficients $
C_{lk}^i$ in terms of a 1-form basis. Equations (\ref{anhol:const}) say that
any vector basis $F_i^\alpha (x)=f_i^\alpha (x)$ such that
its corresponding dual 
1-form basis $\mbox{\boldmath${f^{-1}}^i$}={f^{-1}}^i_\alpha{\rm d}x^\alpha$ 
satisfies the
Maurer-Cartan equation 
\begin{equation}
\label{M-C}{\rm d}\mbox{\boldmath${f^{-1}}^{i}$}=-\frac 12C_{jk}^i
\mbox{\boldmath${f^{-1}}^{j}$}\wedge \mbox{\boldmath${f^{-1}}^{k}$} 
\end{equation}
with constant anholonomicity coefficients $C_{lk}^i={\rm const}$, is a
solution of (\ref{diffWW:1}). The class of all such bases constitutes the
general solution of equations (\ref{diffWW:1}) and therefore the
general solution for 
the coordination bivector satisfying (\ref{diffW:1}) 
and (\ref{AequalW}). $\quad \Box $

The next theorem proves the existence of solutions for the equation (\ref
{diffW:2}) within the class of bivectors 
satisfying (\ref{W}) \cite{Zala:unpub} (its version for a
particular subclass of (\ref{W}), see Section 8 and Proposition
1 below, has been given in \cite{Zala:1992},\cite{Zala:1993}, and
\cite{AZS:1990b} for the case
of Riemannian manifolds).

\begin{theorem}
\label{soldiffW:2} In an arbitrary $n$-dimensional differentiable metric
manifold $({\cal M}$, $g_{\alpha \beta })$ with a volume $n$-form 
there always exist locally volume-preserving bivectors 
${\cal W}_\beta ^{\alpha^{\prime }}(x^{\prime },x)$  
of the form (\ref{W}) with (\ref{C:const}) satisfying (\ref{diffW:2}).
\end{theorem}

{\it Proof.} For a bivector ${\cal W}_\beta ^{\alpha ^{\prime
}}(x^{\prime },x)=f_i^{\alpha ^{\prime }}(x^{\prime })f{^{-1}}_\beta ^i(x)$
equations (\ref{diffW:2}) read 
\begin{equation}
\label{divf}f_{i:\alpha ^{\prime }}^{\alpha ^{\prime }}\equiv f_{i,\alpha
^{\prime }}^{\alpha ^{\prime }}+(\ln \sqrt{-g})_{,\alpha ^{\prime
}}f_i^{\alpha ^{\prime }}=0\ .
\end{equation}
Writing (\ref{divf}) as 
\begin{equation}
-f_i^{\alpha ^{\prime }}f{^{-1}}_{\beta ^{\prime },\alpha ^{\prime }}^i+(\ln 
\sqrt{-g})_{,\beta ^{\prime }}=0\ ,
\end{equation}
and using the definition of the anholonomicity coefficients in terms of 1-form
basis (see (\ref{anhol:const}) above) this equation can be rewritten as:
\begin{equation}
\label{divf2}-f_i^{\alpha ^{\prime }}f{^{-1}}_{\alpha ^{\prime },\beta
^{\prime }}^i-C_{kj}^jf{^{-1}}_{\beta ^{\prime }}^k+(\ln \sqrt{-g})_{,\beta
^{\prime }}=0\ .
\end{equation}
Now it is sufficient to show the integrability of this quasi-linear partial
differential
equation to prove the theorem. Taking into account (\ref
{C:const}) and the property $C_{ji}^iC_{kl}^j=0$ which holds in that case,
one can easily show that (\ref{divf2}) is always integrable on an arbitrary
$n$-dimensional differentiable manifold with a volume 
$n$-form. $\quad \Box $

It should be stressed here that the conditions $C_{jk}^i={\rm const}$ due
to Theorem \ref{soldiffW:1} are essential for the proof of Theorem \ref
{soldiffW:2} and they guarantee the local existence of $n$ linear
independent divergence free vectors, the result holding for both orientable
and non-orientable manifolds (by using the so-called odd volume $n$-form \cite
{deRh:1960},\cite{Mose:1965})\footnote{
$n$ linear independent divergence free vector fields 
on an $n$-dimensional
differentiable manifold with a volume $n$-form with structure functions
(\ref{commff}) can be shown to exist locally
iff $C_{lk,i}^i=0$. This condition is  fulfilled also globally on
parallelizable manifolds (both orientable and non-orientable) \cite
{Grom:1986}. This shows that Theorem \ref{soldiffW:2} is valid globally for
such manifolds (a manifold is called parallelizable if its tangent bundle is
trivial).}.

Another important point is that $n$
vector fields $f_i^\alpha $ satisfying (\ref{commff}) and (\ref{C:const}) 
define a finite dimensional Lie group on the
averaged manifolds which is directly related with the symmetries of such
manifolds \cite{Zala:progr}.

These two theorems prove the existence of solutions for the set of 
equations (\ref{coincidence}), (\ref{idempotency}), (\ref{diffW:1}), 
(\ref{diffW:2}), and (\ref{AequalW}), and therefore prove the existence of the
bilocal operator ${\cal W}_\beta ^{\alpha ^{\prime }}(x^{\prime },x)$ with
the corresponding algebraic and differential properties (\ref{averlocal})
and (\ref{comm}) for the averages (\ref{defaver:MG}).

Now a particular subclass of the operators (\ref{W}) with 
(\ref{C:const}) will be considered to analyze some additional properties 
and to reveal the functional structure of the subclass.

\section{The proper systems of coordinates}

As it has been emphasized in the Introduction and Section 2,
the space-time
averaging procedure adopted in electrodynamics is essentially formulated in
Cartesian coordinates and all its properties are shown by
exploiting the exceptional character of this coordinates (see Section 2 for
details). The covariant formalism developed for the averages (\ref
{defaver:MG}) in macroscopic gravity (Sections 3-6)
generalizes the averaging scheme of
macroscopic electrodynamics for arbitrary $n$-dimensional differentiable 
manifolds and while keeping covariant properties which are analogous
to those in electrodynamics.

Let us now consider the macroscopic gravity averaging scheme for a
particular subclass of operators (\ref{W}) with (\ref{C:const}). 
This particular subclass admits a special coordinate system in which the
averages
and their properties have especially simple form and meaning. Such a
coordinate system is an analogue for macroscopic gravity
of the Cartesian coordinates in Minkowski
space-time.

Let us hereby restrict the class of solutions of the equations (\ref
{diffWW:1}) to the subclass satisfying 
\begin{equation}
\label{coorff}\left[ \mbox{\boldmath$f_{i}$},\mbox{\boldmath$f_{j}$}\right]
=0\,\ , 
\end{equation}
that is $C_{ij}^k\equiv 0$. In this case the vector fields $f_i^\alpha $
constitute a coordinate system and there always exist $n$ functionally
independent scalar functions $\phi ^i(x)$ such that the vector and
corresponding dual
1-form bases are of the form 
\begin{equation}
\label{bases:coor}f_i^\alpha (x(\phi^k))
=\frac{\partial x^\alpha }{\partial \phi ^i}
\,\ ,\quad 
f{^{-1}}_\alpha ^i(\phi(x^\mu))
=\frac{\partial \phi ^i}{\partial x^\alpha }\,\
. 
\end{equation}
Thus, the bilocal operator ${\cal W}_\beta ^{\alpha ^{\prime }}(x^{\prime },x)$
becomes  
\begin{equation}
\label{W:coor}{\cal W}_\beta ^{\alpha ^{\prime }}(x^{\prime },x)=\frac{
\partial x^{\alpha ^{\prime }}}{\partial \phi ^i}\frac{\partial \phi ^i}{
\partial x^\beta }\ . 
\end{equation}
Being functionally independent, the set of $n$ functions $\phi ^i(x)$ can be
taken as a system of local coordinates on the manifold ${\cal M}$ \cite
{Zala:1992},\cite{Zala:1993} (\cite{AZS:1990b} for Riemannian manifolds).

\begin{definition}
\label{propercoor} A coordinate system $\{\phi ^i\}$ defined by $n$
scalar functions $\phi ^i=\phi ^i(x)$ in (\ref{bases:coor}) will be
called a proper coordinate system.
\end{definition}

The usefulness of this definition is motivated by the fact that in a proper
coordinate system the bilocal operator ${\cal W}_\beta ^{\alpha ^{\prime
}}(x^{\prime },x)$ takes the simplest possible form  
\begin{equation}
\label{W:delta}{\cal W}_j^i(\phi ^{\prime },\phi )\equiv {\cal W}_\beta
^{\alpha ^{\prime }}(x^{\prime },x)_{\mid x^\alpha =\phi ^i}=
\mbox{\boldmath$\delta$}_\beta ^{\alpha ^{\prime }}\equiv 
\mbox{\boldmath$\delta$}_j^i, 
\end{equation}
where the bilocal Kronecker symbol $\mbox{\boldmath$\delta$}_\beta ^{\alpha
^{\prime }}$ is defined as $\mbox{\boldmath$\delta$}_\beta ^{\alpha ^{\prime
}}=\delta _i^{\alpha ^{\prime }}\delta _\beta ^i$. The definition of the
average (\ref{defaver:MG}) acquires a remarkable simple form 
(closely
resembling the space-time averages of macroscopic electrodynamics (\ref
{defaver:ED})) when written using a proper coordinate system 
\begin{equation}
\label{defaver:MG:prop}\overline{p}_j^i(\phi )=\frac 1{V_{\Sigma
_\phi}}\int_{\Sigma _\phi}p_j^i(\phi ^{\prime })\sqrt{-g(\phi ^{\prime })}d^n\phi
^{\prime }\ . 
\end{equation}

Theorem \ref{soldiffW:2} has proved the existence of volume-preserving
bilocal operators (\ref{W}) and (\ref{C:const}), and therefore
the existence of
solutions of (\ref{diffW:2}) for  operators ${\cal W}_\beta ^{\alpha
^{\prime }}$ of the form (\ref{W:coor})  from this theorem. It is
useful, however, to prove it here independently due to the above mentioned
importance of the proper system of coordinates for the macroscopic gravity
averaging scheme. The following statement holds in this case 
\cite{Zala:1992},\cite{Zala:1993} (\cite{AZS:1990b} for Riemannian manifolds).

\begin{proposition}
\label{soldiffW:2:coor} In an arbitrary $n$-dimensional differentiable metric
manifold $({\cal M}$, $g_{\alpha \beta })$ with a volume $n$-form, 
there always exist locally a set of $n$ scalar functionally independent
functions
$\phi^{i}(x)$ such that the corresponding coordination bivector
(\ref{W:coor})
satisfies condition (\ref{diffW:2}).
\end{proposition}

{\it Proof.} Let us write down equations (\ref{diffW:2})
for the bivector ${\cal W}_\beta ^{\alpha ^{\prime }}$ (\ref{W:coor}) 
\begin{equation}
\label{vol-pres:eq}\frac{\partial ^2x^{\alpha ^{\prime }}}{\partial \phi
^i\partial \phi ^j}\frac{\partial \phi ^j}{\partial x^{\alpha ^{\prime }}}
+\left( {\rm ln}
\sqrt{-g(x^{\prime })}\right) _{,\mu ^{\prime }}\frac{\partial
x^{\mu ^{\prime }}}{\partial \phi ^i}=0\,  
\end{equation}
in terms of the unknowns $\phi ^i(x)$. These equations 
are equivalent to the system of quasi-linear partial differential equations
(compare with (\ref{divf2})) for $f{^{-1}}_\alpha ^i$ (\ref{bases:coor}) 
\begin{equation}
\label{divf2:coor}-f_i^{\alpha ^{\prime }}(x^{\prime }(\phi^k))
f{^{-1}}_{\alpha^{\prime },j}^i(\phi^k)
+\left({\rm ln} \sqrt{-g(x^{\prime }(\phi ^k))}\right)
_{,j}=0 
\end{equation}
which are always locally integrable on an arbitrary $n$-dimensional 
differentiable manifold
with a volume $n$-form. Then, in accordance with (\ref{bases:coor})
$\phi $\thinspace $^i(x)$ are solutions of the
equations $\phi $\thinspace $_{,\alpha }^i=f{^{-1}}_\alpha ^i$, 
which are also always locally integrable because the integrability
conditions $f{^{-1}}_{[\alpha ,\beta ]}^i=0$ are fulfilled due to 
$C_{jk}^i=0 $ (\ref{coorff}). 
$\quad \Box $

The following result is an obvious consequence of this Proposition. 

\begin{corollary}
\label{corsoldiffW:2:coor}
Any proper coordinate coordinate system such that the corresponding
bivector (\ref{W:coor})
satisfies condition (\ref{diffW:2}) is a volume-preserving system
of coordinates, 
\begin{equation}
\label{vol-pres}
\left( {\rm ln} \sqrt{-g(\phi ^k)}\right)_{,j}=0 \ , \quad 
\mbox{or}, \quad g(\phi ^k)={\rm const} \ . 
\end{equation}
\end{corollary}

It should be noted here that all arguments concerning non-orientable
manifolds and global existence given after Theorem \ref{soldiffW:2} apply
here as well (see Section 7).

For the case of (pseudo)-Riemannian manifolds, Corollary \ref
{corsoldiffW:2:coor} states that in a proper coordinate system the Christoffel
symbols $\Gamma _{\beta \alpha }^\alpha $, (which are $\Gamma _{\beta
\alpha }^\alpha =({\rm ln} \sqrt{-g})_{,\beta }$ due to equi-affinity of Riemannian
manifolds) vanish

\begin{equation}
\label{vol-pres:Riem}\Gamma _i(\phi )\equiv \Gamma _{ij}^j=\left( {\rm ln}
\sqrt{-g(\phi ^k)}\right) _{,j}=0\ . 
\end{equation}

Another useful characterization of the volume-preserving
coordinates, in addition to 
(\ref{vol-pres}), or (\ref{vol-pres:Riem}), can be
obtained in terms of the expansion of the vector fields tangent to the
coordinates lines. Defining a vector $\mbox{\boldmath$\chi_{(i)}$}$ tangent
to a coordinate line $\phi ^i$ as 
\begin{equation}
\label{tangent}\mbox{\boldmath$\chi_{(i)}$}=\chi _{(i)}^j\frac \partial
{\partial \phi ^j}=\delta _{(i)}^j\frac \partial {\partial \phi ^j}\ , 
\end{equation}
it is immediate to find that in the proper coordinate system $\{\phi ^i\}$
the expansion ${\rm div}{\chi }_{(i)}$ of the vector field (\ref{tangent})
is 
\begin{equation}
{\chi }_{(i):j}^j=\left( {\rm ln}\sqrt{-g(\phi ^k)}\right) _{,i}=0\ , 
\end{equation}
so that the condition of vanishing expansion for the tangent vector fields $
\mbox{\boldmath$\chi_{(i)}$}$ is equivalent to the definition of
volume-preserving coordinates (\ref{vol-pres}).

In addition to very simple and transparent forms 
of the coordination bivector (\ref{W:coor}) and the averages
(\ref{defaver:MG:prop}), the volume-preserving proper
coordinate systems also allow a remarkably
simple expression for the commutation between partial differentiation and
averaging. Indeed, expression (\ref{comm}) becomes in the proper
coordinate system
\begin{equation}
\label{commu:2}
\frac \partial {\partial \phi ^k} {\bar{p}_{j}^i(\phi)} =\langle 
\frac \partial {\partial \phi^{\prime k}} {p_{j}^{i}(\phi^\prime)} \rangle \ ,
\end{equation}
which is exactly the same commutation formula for as in 
the averaging scheme in Minkowski manifolds (\ref{comm:ED}).

Let us now study the functional structure of the class of the 
volume-preserving coordinates to understand how large it is and how much 
freedom for coordinate transformations it contains. The following Proposition
reveals the structure of the class.

\begin{proposition}
\label{vp:class} The class of volume-preserving coordinate transformations
on an arbitrary $n$-dimensional differentiable metric manifold 
$({\cal M}$, $g_{\alpha \beta })$ 
with a volume $n$-form,
contains $(n-1)$ arbitrary functions of $n$ arguments
and one arbitrary function of $(n-1)$ arguments.
\end{proposition}

{\it Proof.} Due to Proposition \ref{soldiffW:2:coor} the
volume-preserving coordinates are characterized by the condition (\ref
{vol-pres}). The class of coordinate
transformations $y^\alpha =$ $y^\alpha (x^\mu )$ which preserve this
condition,\thinspace 
\begin{equation}
\label{vp:class:def}\frac \partial {\partial y^\alpha }
\left({\rm ln} \sqrt{-g(y)}\right)=0\ , \quad 
\frac \partial {\partial x^\mu }\left({\rm ln} \sqrt{-g(x)}\right)=0\ , 
\end{equation}
is defined by the system of partial differential equations \cite{ZTE:1996} 
\begin{equation}
\label{vp:class:eqs}\frac{\partial ^2y^\alpha }{\partial x^\nu \partial
x^\mu }\frac{\partial x^\mu }{\partial y^\alpha }=0\ . 
\end{equation}
Equations (\ref{vp:class:eqs}) can be easily written as equations for the
Jacobian of the coordinate transformation $y^\alpha =y^\alpha (x^\mu )$, 
\begin{equation}
\label{vp:class:eqs2}\frac \partial {\partial x^\nu }\left[ \det \left( 
\frac{\partial y^\alpha }{\partial x^\mu }\right) \right] =0\ , 
\end{equation}
which have the general solution 
\begin{equation}
\label{jacob:const}\det \left( \frac{\partial y^\alpha }{\partial x^\mu }
\right) =C 
\end{equation}
where $C$ is a non-vanishing constant to ensure invertibility of the
coordinate transformation. Expanding  the Jacobian (\ref{jacob:const})
by its first row we get 
\begin{equation}
\label{expand}\frac{\partial y^1}{\partial x^1}\xi ^1(y^2,...,y^n)+...+\frac{
\partial y^1}{\partial x^n}\xi ^n(y^2,...,y^n)=C\ 
\end{equation}
where $\xi ^\nu (y^2,...,y^n)$ are the corresponding minors of the
determinant, depending therefore only on $(n-1)$ functions $y^2,...,y^n$.
These  $(n-1)$ functions are 
arbitrary functionally independent functions of $x^\mu $, and $\xi ^\nu 
$ can be taken as explicit functions 
of $x^\mu $, $\xi ^\nu =\xi ^\nu (x^\mu
)$. Then equation (\ref{expand}) reads 
\begin{equation}
\label{expand2}\xi ^\nu (x^\mu )\frac{\partial y^1}{\partial x^\nu }=C\ . 
\end{equation}
This is a linear partial differential equation for the unknown 
$y^1(x^\mu)$. Its general solution is known to be 
\begin{equation}
\label{expand2:sol}y^1(x^\mu )=y_p^1(x^\mu )+h(k^2(x^\mu ),...,k^n(x^\mu )) 
\end{equation}
where $y_p^1(x^\mu )$ is a particular solution of (\ref{expand2}), $
k^A(x^\mu ),$ $A=2,...,n,$ are $n-1$ functionally independent first
integrals satisfying 
\begin{equation}
\label{expand2:integ}\xi ^\nu \frac{\partial k^A}{\partial x^\nu }=0\ , 
\end{equation}
and $h$ is an arbitrary function of $(n-1)$ arguments $k^A$. Thus, the
general solution of the equations (\ref{vp:class:eqs}) defining the class of
volume-preserving coordinates (\ref{vp:class:def}) contains $(n-1)$
arbitrary functions $y^A(x^\mu )$ of $n$ arguments $x^\mu $ and one
arbitrary function $h(k^A)$ of $(n-1)$ variables $k^A(x^\mu )$. $\quad \Box $

The set of proper coordinate systems forms quite a big class and this
functional freedom may be used to specify additional
properties of the averages (\ref{defaver:MG}), or (\ref{defaver:MG:prop}),
when necessary.

Choosing different proper coordinate systems $\phi ^i$ will give different
average fields (\ref{defaver:MG:prop}) of a given microscopic tensor field $
p_\beta ^\alpha (x)$. In general, the averages $\overline{p}_j^i(\phi )$ and 
$\overline{p}_j^i(\widetilde{\phi })$ calculated in the proper coordinate
systems $\phi ^i$ and $\widetilde{\phi }^i$ are not
related by a tensorial law under the transformation $\widetilde{\phi }^i=
\widetilde{\phi }^i(\phi ^j)$, nor are the operators ${\cal W}_j^i(\phi
^{\prime },\phi )$ and ${\cal W}_j^i(\widetilde{\phi }^{\prime },\widetilde{
\phi })$ related by a tensorial transformation with each of them being
$\mbox{\boldmath$\delta$}_j^i$ 
(\ref{W:delta}) in its own proper coordinates system. 
It should be noted here that averages (\ref{defaver:MG}), and (\ref
{defaver:MG:prop}), are obviously tensorial with respect to coordinate
transformations (as follows directly from its definition).
The reason for the ``non-tensorial''
properties between the proper coordinates is due to the structure (\ref{W}),
or (\ref{W:coor}), of the bilocal operator ${\cal W}_\beta ^{\alpha ^{\prime
}}(x^{\prime },x)$ itself, which involves a functional freedom in changing
the functions  $f_i^\alpha (x)$, or $\phi ^i(x)$. This ``non-tensorial''
property is very natural,
indeed, for it states the exceptional character of the proper coordinate
systems for obtaining the simplest and most transparent form
of the averages and the
averaging and coordination operators. It closely resembles, on the other
hand, the definition of averages in macroscopic electrodynamics as it was
noted above, and the exceptional character of Cartesian coordinates used in
that averaging procedure. The class of the proper coordinate systems on an
arbitrary differentiable metric manifold is a natural counterpart of
the Cartesian coordinate system on a Minkowski manifold.
The property they share in common is that both  are volume-preserving.

There is, however, a special subclass within
the class of volume-preserving
coordinate transformations described in Proposition \ref{vp:class} which
keeps the bilocal operator (\ref{W:delta}) and the averages (\ref
{defaver:MG:prop}) covariant.

\begin{proposition}
\label{vp:Lorentz} The class of transformations $\phi ^i\rightarrow
\tilde \phi ^i$ which keeps the bivector ${\cal W}_\beta ^{\alpha ^{\prime }}$
and the averages (\ref{defaver:MG}) covariant within the class of
proper system of coordinates is 
\begin{equation}
\label{Lorentz}\tilde \phi ^i=\Lambda _j^i\phi ^j+a^i
\end{equation}
where $\Lambda _j^i$ and $a^i$ are constant.
\end{proposition}

{\it Proof.} The proof is straightforward. Consider two proper
coordinate systems $\phi ^i$ and $\widetilde{\phi }^i$, and require the
bivectors ${\cal W}_j^i(\phi ^{\prime },\phi )$ and ${\cal W}_j^i(\widetilde{
\phi }^{\prime },\widetilde{\phi })$, each of the form (\ref{W:delta}) in
its own proper coordinates, to be equal 
\begin{equation}
\label{W:delta:Lorentz}\frac{\partial x^{\alpha ^{\prime }}}{\partial \phi ^i
}\frac{\partial \phi ^i}{\partial x^\beta }=\frac{\partial x^{\alpha
^{\prime }}}{\partial \tilde \phi ^i}\frac{\partial \tilde \phi ^i}{\partial
x^\beta }\ . 
\end{equation}
Since the coordinate system $x^\alpha $ in equation (\ref
{W:delta:Lorentz}) is proper by assumption it must be either $\phi ^i$ or $
\widetilde{\phi }^i$. Choosing $\phi ^i$ one immediately obtains
\begin{equation}
\label{W:delta:Lorentz2}\mbox{\boldmath$\delta$}_j^i=\frac{\partial \phi ^i}{
\partial \tilde \phi ^l}_{\mid \tilde \phi (x^{\prime })}\frac{\partial
\tilde \phi ^l}{\partial \phi ^j}_{\mid \phi (x)}\ .\ 
\end{equation}
This last equation has as general solution (\ref{Lorentz}). It is easy to
show that the class of transformations between proper coordinates keeps
the averages (\ref{defaver:MG:prop}) also covariant. $\quad \Box $

Due to Proposition \ref{vp:Lorentz}, if the manifold $({\cal M}$, $g_{\alpha
\beta })$ is chosen to be a (pseudo)-Riemannian spacetime, the averages (\ref
{defaver:MG:prop}) defined in proper coordinates are Lorentz tensors
exactly like the averages in Minkowski space-time (see Section 2).

\section{Conclusion}

Thus, the covariant averaging procedure for objects (tensors, geometric
objects, etc.) in the framework of macroscopic gravity 
\cite{Zala:1996},\cite{Zala:1992},\cite
{Zala:1993} is a natural generalization of the space-time averaging
procedure of macroscopic electrodynamics. It gives a covariant formulation of
the conditions on the averages, which provide them with natural algebraic
and analytical properties. A wide class of averaging and coordination
bilocal operators satisfying all the properties exit locally on an arbitrary
$n$-dimensional differentiable manifold with a volume $n$-form, 
including metric and affine
connection manifolds and, in particular, (pseudo)-Riemannian spaces. The class
of proper coordinate systems, analogous to Cartesian coordinate system 
of Minkowski
space-time and generalizing them,  gives the
simplest and most transparent form of the averages. The averaging procedure
as it is formulated allows a large functional 
freedom which is incorporated
in an elegant way and can be used to arrange
additional specific conditions for the averages.

\section*{Acknowledgments}

MM wishes to thank the
Ministerio de Educaci\'on y Cultura (Spain) for financial support under
grant EX95 40985713. RMZ was supported by a Royal 
Society fellowship (London, UK) and a \'Sniadeccy
Fellowship (Cracow, Poland) and would like to thank the School of
Mathematical Sciences and the N. Copernicus Astronomical Center where the
work has been done in parts for hospitality. RMZ appreciates many
interesting discussions with Boris Dubrovin.

\end{document}